\documentclass[onecolumn]{autart}    

\usepackage[parfill]{parskip}        
\usepackage{graphicx}
\usepackage{amssymb}
\usepackage{epstopdf}
\usepackage[parfill]{parskip}    
\usepackage{amsmath,amssymb,latexsym,enumerate, natbib, mathrsfs}
\usepackage{graphicx}
\usepackage{array}
\usepackage{times}
\usepackage{moreverb}
\usepackage{marvosym}
\usepackage{hyperref}

\DeclareMathOperator*{\argmax}{arg\,max}

\newtheorem{remark}{Remark}

\newtheorem{proposition}{Proposition}

\begin{document}

\begin{frontmatter}

\title{Optimal residence time control for stochastically perturbed prescription opioid epidemic models\thanksref{footnoteinfo}} 

\thanks[footnoteinfo]{Corresponding author Getachew K Befekadu. Tel. +1 443 885 3073.
	          Fax +1 443 885 8219.}

\author[GKB]{Getachew K Befekadu}\ead{getachew.befekadu@morgan.edu}   

\address[GKB]{Department of Electrical \& Computer Engineering, Clarence M. Mitchell, Jr. School of Engineering, Morgan State University, 1700 E. Cold Spring Lane, Schaefer Engr. Bldg. 331, Baltimore, MD 21251, USA.}
%
       
\begin{keyword}                           
Diffusion processes; exit probability;  epidemiology; population compartmental model; prescription drug addiction; Markov controls; minimum exit-rates; optimal control problem.
\end{keyword}                             

\begin{abstract}                          
In this paper, we consider an optimal control problem for a prescription opioid epidemic model that describes the interaction between the regular prescription or addictive use of opioid drugs, and the process of rehabilitation and that of relapsing into opioid drug use. In particular, our interest is in the situation, where the control appearing linearly in the opioid epidemics is interpreted as the rate at which the susceptible individuals are effectively removed from the population due to an opioid-related intervention policy or when the dynamics of the addicted is strategically influenced due to an accessible addiction treatment facility, while a small perturbing noise enters through the dynamics of the susceptible group in the population compartmental model. To this end, we introduce a mathematical apparatus that minimizes the asymptotic exit-rate with which the solution for such stochastically perturbed prescription opioid epidemics exits from a given bounded open domain. Moreover, under certain assumptions, we also provide an admissible optimal Markov control for the corresponding optimal control problem that optimally effected removal of the susceptible or recovered individuals from the population dynamics.
\end{abstract}

\end{frontmatter}

\section{Introduction} \label{S1}
Despite the urgency of the problem, relatively little is known about how to address the current opioid epidemics based on a systematic multi-pronged approach aiming at a wide spectrum issues that arise from misuse of prescription drugs, inappropriate opioid prescribing practices and higher prescribing rates, or due to the lack of easily accessible opioid dependence treatment facilities (e.g., see \cite{r1}, \cite{r2}, \cite{r3}, \cite{DowHC16} and \cite{OASPE15} for general policies that are primarily aimed at curbing prescription opioid abuse, preventing inappropriate prescribing practices, developing abuse deterrents or preventing drug diversion mechanisms). In this paper, without attempting to give a complete literature review, we make an effort to address one aspect of this complex problem using a combination of optimal control theory and epidemiological insights. In particular, we consider an optimal control problem for a prescription opioid epidemic model that describes the interaction between the regular prescription or addictive use of opioid drug, and the process of rehabilitation and that of relapsing into opioid drug use. Our main interest is in the situation, where the control appearing linearly in the opioid epidemic model is interpreted as the rate at which the susceptible individuals are effectively removed from the population due to an opioid-related intervention policy or when the dynamics of addicted is strategically influenced due to an accessible addiction treatment facility, while a small perturbing noise enters through the dynamics of the susceptible group in the population compartmental model (e.g., see \cite{BefZ18} for additional discussions how random perturbing noise may propagate through the dynamics of prescription opioid epidemics). To this end, we introduce a mathematical apparatus that  minimizes the asymptotic exit-rate (equivalently, maximizes  the expected exit-time or the residence time) with which the solution for such stochastically perturbed prescription opioid epidemics exits from a given bounded open domain. Moreover, under certain assumptions, we also provide an admissible optimal Markov control for the corresponding optimal control problem that optimally effected removal of the susceptible or the recovered from the population dynamics. 

Finally, it is worth mentioning that some interesting studies on the minimum escape time or optimal residence time control problem for stochastically perturbed dynamical systems have been reported in literature (e.g., see \cite{DuFre98}, \cite{DuKu98}, \cite{KimMR90} or \cite{BefA15a}, among others, mainly from a mathematical control theory point of view). The rationale behind our framework, which follows in some sense the settings of these papers, is to provide a stochastic optimal control argument which is more practical for characterizing typical sample paths of regular prescription opioid users, opioid addicts or the process of rehabilitation and relapsing into opioid drug uses, with some form of opioid-related intervention strategies or policies.

The remainder of this paper is organized as follows. In Section~\ref{S2}, we present our problem formulation for the optimal residence time control of stochastically perturbed prescription opioid epidemic model. In Section~\ref{S3}, we present our main results -- where we provide mathematical arguments that characterize the admissible solutions for the optimal residence time control problem. Section~\ref{S4} contains simulation results. Finally, Section~\ref{S5} provides concluding remarks.

\section{Methods} \label{S2}
In this section, we present our problem formulation, where our interest is to minimize the asymptotic exit-rate (or equivalently maximizing the expected exit-time) with which the solution of stochastically perturbed prescription opioid epidemic model exits from the given bounded open domain.

\subsection{Model description} \label{S2(1)}
In what follows, if we denote the susceptible, addicted and recovered in a closed population by $X_1(t)$, $X_2(t)$ and $X_3(t)$, respectively. Then, the prescription opioid epidemic dynamical model, with small random perturbing noise, can be written as follows\footnote{Note that, for a normalized population, i.e., when $N=1$, the proportion of opioid prescription users, denoted by $Z(t)$, is given by $Z(t) = 1 - X_1(t)-X_2(t)-X_3(t)$, for $t \ge 0$.}
\begin{eqnarray}
\left.\begin{array}{l}
dX_1(t) = f_1(X_1(t), X_2(t), X_3(t)) dt + \sqrt{\epsilon} dW(t)\\
dX_2(t) = f_2(X_1(t), X_2(t), X_3(t)) dt \\
dX_3(t) = f_3(X_2(t), X_3(t)) dt
\end{array}\right\} \label{Eq2.6}
\end{eqnarray}
where $\bigl(W(t)\bigr)_{t \ge 0}$ is a one-dimensional Brownian motion, $\bigl(X_1(t), X_2(t), X_3(t)\bigr)_{t \ge 0}$ being an $\mathbb{R}^3$-valued degenerate diffusion process, and
$\epsilon$ is a small positive number that represents the level of the random perturbation in the prescription opioid epidemic dynamics, while the functions $f_1$, $f_2$ and $f_3$ are given by
\begin{align*}
f_1(t, x_1, x_2, x_3)  &=  -\alpha x_1(t) -\beta(1-\xi)x_1(t) x_2(t) - \beta \xi x_1(t) (1 - x_1(t) - x_2(t) - x_3(t)) \\
                                & \quad\quad  + (\varepsilon + \mu)(1 - x_1(t) - x_2(t) - x_3(t)) + (\delta + \mu) x_3(t) + \mu^{\ast} x_2(t),
\end{align*}
\begin{align*}
f_2(t, x_1, x_2, x_3) &= \gamma (1 - x_1(t) - x_2(t) - x_3(t)) + \sigma x_3(t) + \beta(1 - \xi) x_1(t) x_2(t)  \\
                                 & \quad \quad  + \beta \xi x_1(t) (1 - x_1(t) - x_2(t) - x_3(t)) + \nu x_3(t) x_2(t) - (\zeta + \mu^{\ast}) x_2(t)
\end{align*}
and
\begin{align*}
f_3(t, x_2, x_3) = \zeta x_2(t) - \mu x_3(t) x_2(t) - (\delta + \sigma + \mu) x_3(t),
\end{align*}
respectively (e.g., see also \cite{r6} and \cite{BefZ18} for additional discussions on the detailed model derivation). In the above prescription opioid epidemic dynamical population model, we assume that no new addictive opioid drug users are introduced from outside, but there is an external small random perturbing noise that enters through the dynamics of the susceptible group and then its effect is subsequently propagated to the other subsystems. Moreover, Table~\ref{TB1} contains a brief description of the system parameters in Equation~\eqref{Eq2.6}.

Furthermore, if we denote by a bold letter a quantity in $\mathbb{R}^3$, for example, the solution in Equation~\eqref{Eq2.6} is denoted by $\bigl(\mathbf{X}(t)\bigr)_{t \ge 0} = \bigl(X_1(t), X_2(t), X_3(t)\bigr)_{t\ge 0}$, then we can rewrite Equation~\eqref{Eq2.6} as follows
\begin{align}
d \mathbf{X}(t) = \mathbf{F} (\mathbf{X}(t)) dt + \sqrt{\epsilon} B dW(t), \label{Eq2.7}
\end{align}
where $\mathbf{F} = \bigl[f_1, f_2, f_3\bigr]^T$ is an $\mathbb{R}^3$-valued function and $B$ stands for a column vector that embeds $\mathbb{R}$ into $\mathbb{R}^3$, i.e., $B = [1, 0, 0]^T$. Note that the corresponding degenerate elliptic operator for the diffusion process $\mathbf{X}(t)$ is given by
\begin{align}
\mathcal{L}^{\epsilon} (\cdot) (\mathbf{x}) = \frac{\epsilon}{2} \operatorname{tr}\Bigl \{D_{x_1}^2 (\cdot) \Bigr\} + \sum\nolimits_{i=1}^3 f_{i}(\mathbf{x}) D_{x_i} (\cdot),  \label{Eq2.8}
\end{align}
where $D_{x_i}$ and $D_{x_1}^2$ (with $D_{x_1}^2 = \bigl({\partial^2 }/{\partial x_1 \partial x_1} \bigr)$) are the gradient and the Hessian, w.r.t. the variable $x_i$, for $i \in \{1,2,3\}$, respectively.

\begin{table} [h!]
\begin{center}
\caption{Summary of notation} \label{TB1}
\begin{tabular}{| p{1.25cm} | p{10cm} |}
\hline
$\alpha$ & the rate at which people are prescribed opioids\\ 
\hline
$\beta$ & the total probability of becoming addicted to opioids other than by prescription\\
\hline
$\beta(1-\xi)$ & the proportion of $\beta$ caused by black market drugs or other addicts \\
\hline
$\beta \xi$ & the rate at which the non-prescribed, susceptible individuals begin abusing opioids due to the accessibility of extra prescription opioids, e.g., new addict users may get drugs from a friend or relative's prescription \\
\hline
$\varepsilon$ & the rate at which people come back to the susceptible group after being prescribed opioids\\
\hline
$\delta$ & the rate at which people come back to the susceptible group after successfully finishing treatment. Despite having completed rehabilitation, we assume people are susceptible to addiction for life \\
\hline
$\mu$ & the natural death rate \\
\hline
$\mu^{\ast}$ & the (enhanced) death rate for addicts ($\mu$ plus overdose rate)\\
\hline
$\gamma$ & the rate at which the prescribed opioid users fall into addiction\\
\hline
$\zeta$ & the rate at which addicted/dependent opioid users enter the treatment/rehabilitation process\\
\hline
$\nu$ & the rate at which users during the treatment fall back into addictive drug use due to the availability of prescribed painkillers from relatives or friends\\
\hline
\end{tabular}
\end{center}
\end{table}

\subsection{Controlled-eigenvalue problem} \label{S2(2)}
In this subsection, we consider the following controlled version of SDE for Equation~\eqref{Eq2.7}, with the corresponding controlled-diffusion process $\bigl(\mathbf{X}_{0,\mathbf{x}}^{u,\epsilon}(t) \bigr)_{t \ge 0}$, i.e., 
\begin{align}
d \mathbf{X}_{0,\mathbf{x}}^{u,\epsilon}(t) = \bigl[\mathbf{F} (\mathbf{X}_{0,\mathbf{x}}^{u,\epsilon}(t)) + \tilde{B} u(t) \bigr]dt +  \sqrt{\epsilon} B dW(t), \,\, \mathbf{X}_{0,\mathbf{x}}^{u,\epsilon}(0) = \mathbf{x}, \label{Eq2.15}
\end{align}
where $u(\cdot)$ is a measurable control process from a set $\mathcal{U}$ which is $\mathbb{R}^2$-valued progressively measurable processes (i.e., a family of nonanticipative processes, for all $t > s$, $(W(t)-W(s))$ is independent of $u(r)$ for $r \le s$) and such that

\begin{align*}
\mathbb{E} \int_{0}^{\infty} \vert u(t)\vert^2 dt < \infty,
\end{align*}
and the matrix $\tilde{B}$ is given by 
\begin{align*}
\tilde{B} = \left[ \begin{array}{cc} b_1 ~&~ 0\\
  0 ~&~ 0 \\
  0 ~&~ b_2 \end{array} \right],
\end{align*}
while the numerical values for $b_j \ge 0$, with $j=1,\,2$, describe the efficiency or effectiveness of the control efforts.\footnote{Note that small $b_1$ and $b_2$ also imply that the admissible control $u \in \mathcal{U}$ is {\it expensive} due to some costs associated with its implementation.} 

\begin{remark}  \label{R1}
Note that, from the structure of matrix $\tilde{B}$, we also observe that the admissible control strategy influences the opioid epidemic dynamics directly through the susceptible or addictive groups, where such an admissible control can be interpreted as the rate at which the susceptible individuals are effectively removed from the population or when the dynamics of addicted is strategically influenced due to an easily accessible addiction treatment facility.
\end{remark}

Let $D \subset \mathbb{R}^3$ be a given bounded open domain, with smooth boundary $\partial D$ (i.e., $\partial D$ is a manifold of class $C^2$), and let us denote by $C^{\infty}(D)$ the spaces of infinitely differentiable functions on $D$. Furthermore, let $\mathbb{P}_{\mathbf{x}}^{u,\epsilon} \bigl\{\mathcal{A}\bigr \}$ and $\mathbb{E}_{\mathbf{x}}^{u,\epsilon}\bigl\{\omega \bigr\}$, as usual, denote the probability of an event $\mathcal{A}$ and the expectation of a random variable $\omega$, respectively, for the controlled-diffusion process $\mathbf{X}_{0,\mathbf{x}}^{u,\epsilon}(t)$ starting from $\mathbf{x} \in D$. 

In this paper, our main interest is to confine the controlled-diffusion process $\mathbf{X}_{0,\mathbf{x}}^{u,\epsilon}(t)$ inside the domain $D$ as long as possible. A standard formulation for such a problem is to maximize the expected exit-time (or the residence time) with respect to a certain class of admissible controls, i.e.,
 \begin{align*}
 \max_{u \in \mathcal{U}}\mathbb{E}_{\mathbf{x}}^{u,\epsilon}\bigl\{\tau_{D}^{\epsilon}\bigr\}, 
 \end{align*}
where $\tau_{D}^{\epsilon}$ denotes the first exit-time for the controlled-diffusion process $\mathbf{X}_{0,\mathbf{x}}^{u,\epsilon}(t)$ from the domain $D$ and is given by 
 \begin{align}
\tau_{D}^{\epsilon} = \inf \Bigl\{ t > 0 \, \bigl\vert \, \mathbf{X}_{0,\mathbf{x}}^{u,\epsilon}(t) \in \partial D \Bigr\}. \label{Eq-2}
\end{align}
Note that, in general, it is difficult to get effective information about the exit probability (or the maximum expected exit-time) and, at the same time, a set of admissible controls in this way. On the other hand, one could consider a more natural objective that is directly associated with the asymptotic exit-rate with which the controlled-diffusion process $\mathbf{X}_{0,\mathbf{x}}^{u,\epsilon}(t)$ exits from the domain $D$. Furthermore, this suggests minimizing the principal eigenvalue $\lambda_{v}^{\epsilon}$ 
\begin{align}
\lambda_{u}^{\epsilon} = - \limsup_{t \rightarrow \infty} \frac{1} {t} \log \mathbb{P}_{\mathbf{x}}^{u,\epsilon} \bigl\{\tau_{D}^{\epsilon} > t \bigr\}, \label{Eq1.6}
\end{align}
with respect to a certain class of admissible controls.\footnote{Recently, the authors in \cite{AraBB16} have provided some interesting results on controlled equilibrium selection problem in stochastically perturbed dynamics, but in a different context.} 

In what follows, we specifically consider a {\em precise} stationary Markov control $u(t) = v\bigl(\mathbf{X}_{0,\mathbf{x}}^{v,\epsilon}(t)\bigr) \in \mathcal{U}$, for $t \ge 0$, with some measurable map $v \colon \mathbb{R}^{3} \rightarrow \mathcal{U}$. Then, we suppose that the controlled-SDE in Equation~\eqref{Eq2.15} is composed with this admissible Markov control $v$. Note that the extended generator for the controlled-diffusion process $\mathbf{X}_{0,\mathbf{x}}^{v,\epsilon}(t)$ is given by
\begin{eqnarray}
\mathcal{L}_{v}^{\epsilon} \bigl(\cdot\bigr) \bigl(\mathbf{x}\bigr) = \frac{\epsilon}{2} \operatorname{tr}\Bigl \{D_{x_1}^2 (\cdot) \Bigr\} + \Bigl\langle \mathbf{F}(\mathbf{x}) + \tilde{B} v(\mathbf{x}), D_{\mathbf{x}}(\cdot) \Bigr\rangle, \label{Eq-4}
\end{eqnarray}
where $D_{\mathbf{x}}(\cdot)$ denotes the gradient operator, i.e., $D_{\mathbf{x}}(\cdot) \equiv [D_{x_1}(\cdot),\,D_{x_2}(\cdot),\,D_{x_3}(\cdot)]^T$, with respect to $\mathbf{x} \in \mathbb{R}^3$. 

Next, let us consider the following controlled-eigenvalue problem
\begin{align}
\left.\begin{array}{c}
  - \mathcal{L}_{v}^{\epsilon} \psi_{v} \bigl(\mathbf{x}\bigr) = \lambda_{v}^{\epsilon} \psi_{v}\bigl(\mathbf{x}\bigr) \quad \text{in} \quad D \quad \vspace{2mm} \\
  \quad \psi_{v}\bigl(\mathbf{x}\bigr) = 0 \quad \text{on} \quad \partial D  \label{Eq-7}
  \end{array}\right\} 
\end{align}
where the extended generator $\mathcal{L}_{v}^{\epsilon}$ is given in Equation~\eqref{Eq-4} above. 

In the following section, using Theorems~1.1, 1.2 and 1.4 from \cite{QuaSi08} (see also \cite[Proposition~3.2]{BefA15a}), we provide a condition for the existence of a unique principal eigenvalue $\lambda_{v}^{\epsilon} > 0$ and an eigenfunction $\psi_{v} \in W_{loc}^{2,p} \bigl(D\bigr) \cap C\bigl(\bar{D}\bigr)$ pairs for the eigenvalue problem in Equation~\eqref{Eq-7}, with zero boundary condition on $\partial D$. Then, we further make use of the following observation -- where such an eigenvalue $\lambda_{v}^{\epsilon}$ is also related to the minimum asymptotic exit-rate with which the controlled-diffusion process $\mathbf{X}_{0,\mathbf{x}}^{v,\epsilon}(t)$ exits from the given domain $D$, when the controlled-SDE in Equation~\eqref{Eq2.15} is composed with an admissible Markov control $v$.

\section{Characterizing the optimal residence time control problem} \label{S3}
In this section, we present our main result that characterizes the admissible solutions for the optimal control problem in Equation~\eqref{Eq-7}. In particular, Proposition~\ref{P-1} establishes a connection between the minimum asymptotic exit-rate (i.e., the optimal residence time) and the principal eigenvalue problem for the extended generator $\mathcal{L}_{v}^{\epsilon}$ in Equation~\eqref{Eq-4}. Whilst, in Remark~\ref{R2} provides a condition for the existence of an admissible optimal Markov control for the corresponding optimal control problem.

\begin{proposition} \label{P-1}
Suppose that an admissible Markov control $v$ is given, then the principal eigenvalue $\lambda_{v}^{\epsilon}$ for the extended generator $\mathcal{L}_{v}^{\epsilon}$, with zero boundary condition on $\partial D$, is given by 
\begin{align}
\lambda_{v}^{\epsilon} = - \limsup_{t \rightarrow \infty} \frac{1} {t} \log \mathbb{P}_{\mathbf{x}}^{v,\epsilon} \bigl\{\tau_{D}^{\epsilon} > t \bigr\},  \label{Eq-8}
\end{align}
where $\tau_{D}^{\epsilon}$ is the first exit-time for the controlled-diffusion process $\mathbf{X}_{0,\mathbf{x}}^{v,\epsilon}(t)$ from the given bounded domain $D$, while the probability $\mathbb{P}_{\mathbf{x}}^{v,\epsilon} \bigl\{\cdot\bigr\}$ in Equation~\eqref{Eq-8} is conditioned on the initial condition $\mathbf{x} \in D$ as well as on the admissible Markov control $v(\cdot)$ for $t \in [0,\, \tau_{D}^{\epsilon})$.
\end{proposition}

\begin{remark} \label{R1}
Note that the above proposition establishes a connection between the minimum exit rate problems and that of the principal eigenvalue problems, where such a connection is also studied in \cite{BisB10}.
\end{remark}

In what follows, let us define the following Hamilton-Jacobi-Bellman (HJB) equation
\begin{align}
&\mathcal{L}_{u}^{\epsilon} \bigl(\cdot\bigr) \bigl(\mathbf{x}, u\bigr) = \frac{\epsilon}{2} \operatorname{tr}\Bigl \{D_{x_1}^2 (\cdot) \Bigr\} + \Bigl\langle \mathbf{F}(\mathbf{x}) + \tilde{B}u, D_{\mathbf{x}}(\cdot) \Bigr\rangle, \label{Eq-10}
\end{align}
with $D_{\mathbf{x}}(\cdot) = [D_{x_1}(\cdot),\,D_{x_2}(\cdot),\,D_{x_3}(\cdot)]^T$.

\begin{remark} \label{R2}
 Here, it is worth remarking that one can also associate the above HJB equation with the following optimal control problem
\begin{align}
\max_{u \in \mathbb{R}} \Bigl \{ \mathcal{L}_{u}^{\epsilon}\psi\bigl(\mathbf{x}, u\bigr) + \lambda \psi\bigl(\mathbf{x}\bigr) \Bigr\}, \label{Eq-11}
\end{align}
where the admissible optimal control $u^{\ast}()$ can be determined by any measurable selector of 
\begin{align*}
\argmax \bigl\{\mathcal{L}_{u}^{\epsilon} \psi \bigl(\mathbf{x}, \,\cdot \,\bigr) \bigr\}, \quad \mathbf{x} \in D.
\end{align*}
\end{remark}

\section{Simulation results} \label{S4}
In this section, we apply our framework that is discussed in the previous section. In particular, we determined an optimal Markov control strategy for the prescription opioid epidemics with small random perturbing noise. Note that the control strategy appears linearly in the model and, as a result of this, we interpreted such an admissible control strategy as the rate at which the susceptible individuals are effectively removed from the population due to an intervention strategy or when the dynamics of addicted is influenced due to an accessible addiction treatment facility. In the simulation results, we used literature based parameter values (that are given in Table~\ref{TB2}) for the prescription opioid epidemic model (e.g., see also \cite{r6} for additional discussions). Here, we are mainly interested in the addiction-free equilibrium case, i.e., for $\gamma = 0$, $\xi = 0$ and $\beta > 0$, when the linearized prescription opioid epidemic model (corresponding to the deterministic model, i.e., $\dot{\tilde{\mathbf{X}}}(t) = \mathbf{F} (\tilde{\mathbf{X}}(t))$) becomes an addiction-free equilibrium, with the following steady state conditions\footnote{The number of opioid prescription users corresponding to an addiction-free equilibrium point is computed as follows $Z^{\ast} = 1 - X_1^{\ast} - X_2^{\ast} - X_3^{\ast}$, that is, $Z^{\ast} = 0.0565$, while the addiction-free equilibrium is given by $(X_1^{\ast}, X_2^{\ast}, X_3^{\ast}) = (0.9435,\,0,\, 0)$.}
\begin{align*}
X_1^{\ast} = \frac{\varepsilon + \mu} {\alpha + \varepsilon + \mu}, ~\quad X_2^{\ast} = 0, ~\quad X_3^{\ast} = 0 ~~ \text{and } ~~ Z^{\ast} =\frac{\alpha} {\alpha + \varepsilon + \mu}.
\end{align*}

\begin{table} [h!]
\begin{center}
\caption{Literature based parameter values(see \cite{r6})} \label{TB2}
\begin{tabular}{ c c | c c}
\hline\hline
\rm Parameter & \rm Numerical value & \rm Parameter  & \rm Numerical value \\ 
\hline
$\alpha$ & \rm 0.15 &  $\delta$ & \rm 0.1\\ 
$\varepsilon$ & \rm 0.8 - 8 & $\nu$  & \rm 0.2\\ 
$\beta$ & \rm 0.0036 &  $\sigma$  & \rm 0.7\\ 
$\xi$ & \rm 0.74 & $\mu$  & \rm 0.007288\\ 
$\gamma$ & \rm 0.00744 & $\mu^{\ast}$  & \rm 0.01155\\ 
$\zeta$ & \rm 0.2 - 2 &  \rm -  & \rm -\\ 
\hline
\end{tabular}
\end{center}
\end{table}

Note that the Jacobian matrix $J(\mathbf{X})$, i.e., the linearized prescription opioid epidemic model when evaluated at the addiction-free equilibrium point $\mathbf{X}^{\ast}$, is given by
\begin{align*}
  J(\mathbf{X}) \bigl \vert_{\mathbf{X}=\mathbf{X}^{\ast}} &= \left[\frac{\partial f_i(\mathbf{X})}{\partial X_j}\right]_{ij} \biggl \vert_{\mathbf{X} =\mathbf{X}^{\ast}}, \qquad i,j \in \{1,2,3\} \\
                                  &= \left[ \begin{array}{ccc} -(\alpha + \varepsilon + \mu) & ~~~~ \dfrac{\beta(\varepsilon+\mu)}{\alpha + \varepsilon+ \mu} - (\varepsilon + \mu) + \mu^{\ast} ~~~~ & \delta - \varepsilon\\
  0 & \dfrac{\beta(\varepsilon+\mu)}{\alpha + \varepsilon+ \mu} - (\zeta + \mu^{\ast}) & \sigma \\
  0 & \zeta  & -(\delta + \sigma + \mu)
  \end{array} \right]
\end{align*}
and the corresponding eigenvalues for the Jacobian matrix $J(\mathbf{X^{\ast}})$, that is, $\big\{-3.1573, -0.0323, -1.0331\big\}$, are all strictly negative and, hence, the addiction-free equilibrium is asymptotically stable, with a reproduction number $\mathcal{R}_{\rm o} = 0.0766.$\footnote{Note that, if the reproduction number $\mathcal{R}_{\rm o} > 1$ or $\mathcal{R}_{\rm o} < 1$, then it evidently indicates the prevalence of opioid addicts in the population.} Moreover, we performed our simulation studies using parameter values from Table~\ref{TB2}, with $\varepsilon=3$ and $\zeta=0.25$. Note that, from the physical point of view, the domain of interest, i.e., a bounded open domain $D \subset \mathbb{R}^3$, with smooth boundary condition, must be contained inside the following boundary condition 
\begin{align*}
D \subset \left\{\begin{array}{c}
X_i(t) \ge 0, \quad\quad i=1,2,3\\
X_1(t) + X_2(t) + X_3(t) +Z(t) = 1, \quad \forall t \ge 0
\end{array}\right\} 
\end{align*}
or equivalently, with $Z(t) \ge 0$, for all $t \ge 0$,
\begin{align*}
D \subset \left\{\begin{array}{c}
X_i(t) \ge 0, \quad\quad i=1,2,3\\
X_1(t) + X_2(t) + X_3(t) \le 1, \quad \forall t \ge 0
\end{array}\right\},
\end{align*}
with smooth boundary $\partial D$.

In what follows, we provided an upper bound for the optimal residence time based on the above linearized prescription opioid epidemic model.\footnote{Here, we assume that the dynamics of the linearized system is close to that of the original prescription opioid epidemic model in Equation~\eqref{Eq2.6}.} In the simulation, we set $b_1=0.01$ and $b_2=0.001$ to highlight the relative effectiveness of those parameters in the controller matrix $\tilde{B}$. Note that, for the control strategy appearing linearly in the opioid epidemic model, the controllability property of the system holds true (see \cite{BefZ18} for additional discussions for such an assumption). As a result of this, there exists at least one Markovian control $v^{\ast}$, with eigenvalue-eigenfunction pair $(\lambda_{v^{\ast}}, \varphi^{\ast})$ (cf. Proposition~\ref{P-2} above). Furthermore, for any small noise intensity $\epsilon \ll 1$,
the corresponding first exit-time $\tau_D^{\epsilon}$ for the controlled-diffusion process $(\mathbf{X}_{0,\mathbf{x}}^{v^{\ast},\epsilon}(t))_{t \ge 0}$ with which it exits from the domain $D$ is finite and always bounded from the above, i.e., 
\begin{align*}
\tau_D^{\epsilon} \le \tilde{\tau}_D^{\epsilon}\bigl(\tilde{\mathbf{X}}, K_{\tilde{\gamma}}\bigr),
\end{align*}
where $\tilde{\tau}_D^{\epsilon} \bigl(\tilde{\mathbf{X}}, K_{\tilde{\gamma}}\bigr)$ denotes the first exit-time for the diffusion process $\tilde{\mathbf{X}}_{0,\mathbf{x}}^{K_{\tilde{\gamma}},\epsilon}(t)$ from the domain $D$ corresponding to the linearized prescription opioid epidemic model with small perturbing noise, i.e.,
\begin{align*}
\tilde{\tau}_D^{\epsilon}\bigl(\tilde{\mathbf{X}}, K_{\tilde{\gamma}}\bigr) = \inf \bigl\{ t > 0 \, \bigl\vert \, \tilde{\mathbf{X}}_{0,\mathbf{x}}^{K_{\tilde{\gamma}},\epsilon}(t) \in \partial D \bigr\},
\end{align*}
 while the diffusion process $(\tilde{\mathbf{X}}_{0,\mathbf{x}}^{K_{\tilde{\gamma}},\epsilon}(t))_{t \ge 0}$ satisfies the following the controlled-linear SDE
\begin{align*}
 d \tilde{\mathbf{X}}_{0,\mathbf{x}}^{K_{\tilde{\gamma}},\epsilon}(t) = \bigl(\mathbf{A} + \tilde{B} K_{\tilde{\gamma}} \bigr) \tilde{\mathbf{X}}_{0,\mathbf{x}}^{K_{\tilde{\gamma}},\epsilon}(t) dt + \sqrt{\epsilon} B d W (t),
\end{align*}
with $\mathbf{A}=J(\mathbf{X}) \bigl \vert_{\mathbf{X}=\mathbf{X}^{\ast}}$ and $K_{\tilde{\gamma}}$ is some state feedback control matrix that is given below. Moreover, the logarithm residence time satisfies the following
\begin{align*}
 \epsilon \ln \tilde{\tau}_D^{\epsilon}\bigl(\tilde{\mathbf{X}}, K_{\tilde{\gamma}}\bigr) = \tilde{\phi}(D, K_{\tilde{\gamma}}),
\end{align*}
where
\begin{align*}
 \tilde{\phi}(D, K_{\tilde{\gamma}}) = \inf_{\mathbf{X} \in \partial D} \frac{1}{2\Vert \mathbf{P}_{\tilde{\gamma}} \Vert } \mathbf{X}^T \mathbf{P}_{\tilde{\gamma}} \mathbf{X}
\end{align*}
and $\mathbf{P}_{\tilde{\gamma}}$ is a p.d.s. matrix that further satisfies the following algebraic equality equation\footnote{Similar results, based on large deviations principles, for general linear systems, with constant diffusion terms, have been studied in detail by some authors (e.g., see \cite{DuFre98}, \cite{DuKu98} or \cite{KimMR90}).}
\begin{align*}
\mathbf{A}^T \mathbf{P}_{\tilde{\gamma}}  + \mathbf{P}_{\tilde{\gamma}} \mathbf{A} + I - \frac{1}{\tilde{\gamma}} \mathbf{P}_{\tilde{\gamma}} \tilde{B} \tilde{B}^T \mathbf{P}_{\tilde{\gamma}}=0, \quad \tilde{\gamma} \ge 0, 
\end{align*}
while the feedback control matrix is given by $K_{\tilde{\gamma}} = (-1/{\tilde{\gamma}}) \tilde{B}^T \mathbf{P}_{\tilde{\gamma}}$.

\begin{figure}[h!]
\begin{center}
{\includegraphics[scale=0.311]{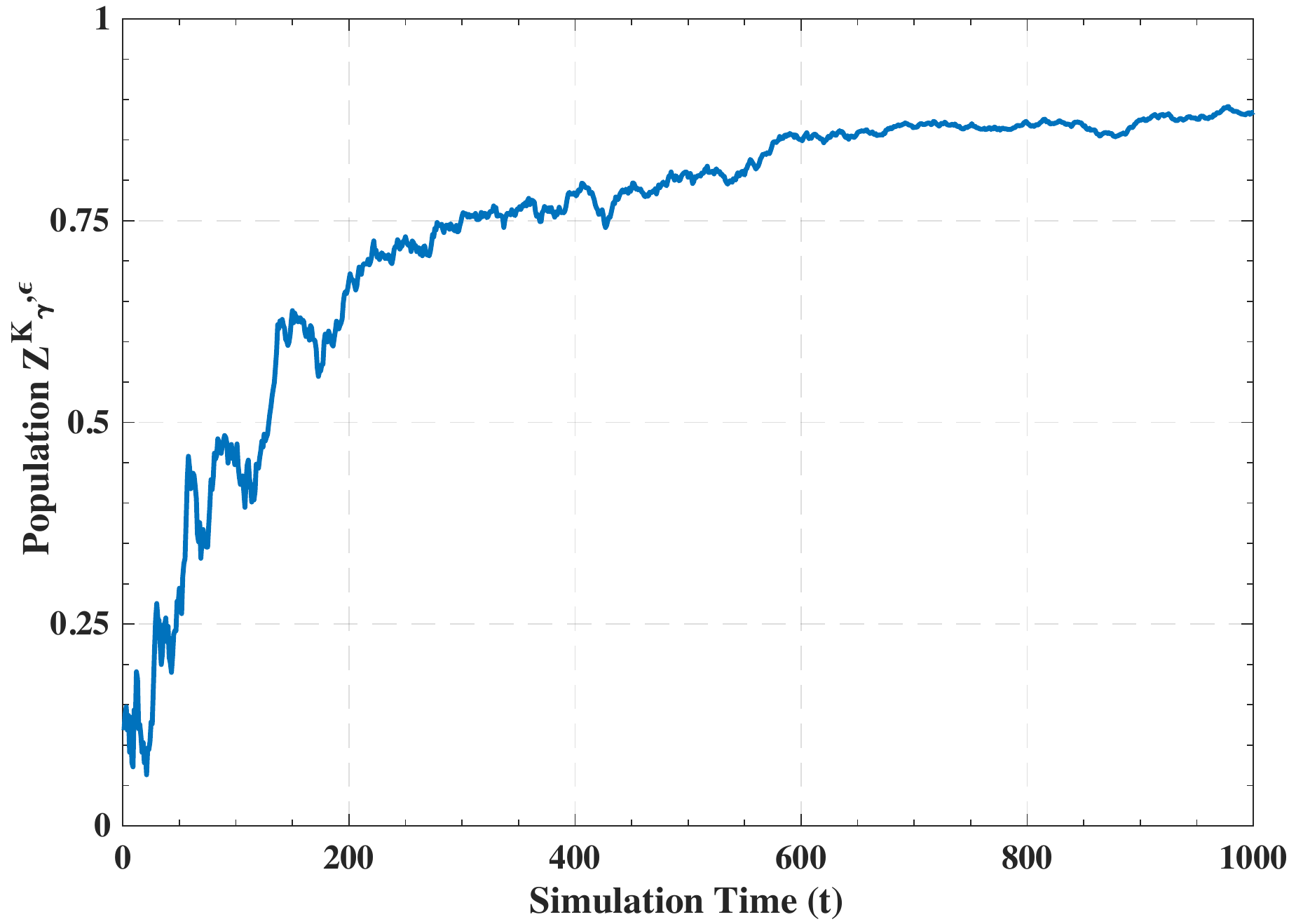}}
{\includegraphics[scale=0.311]{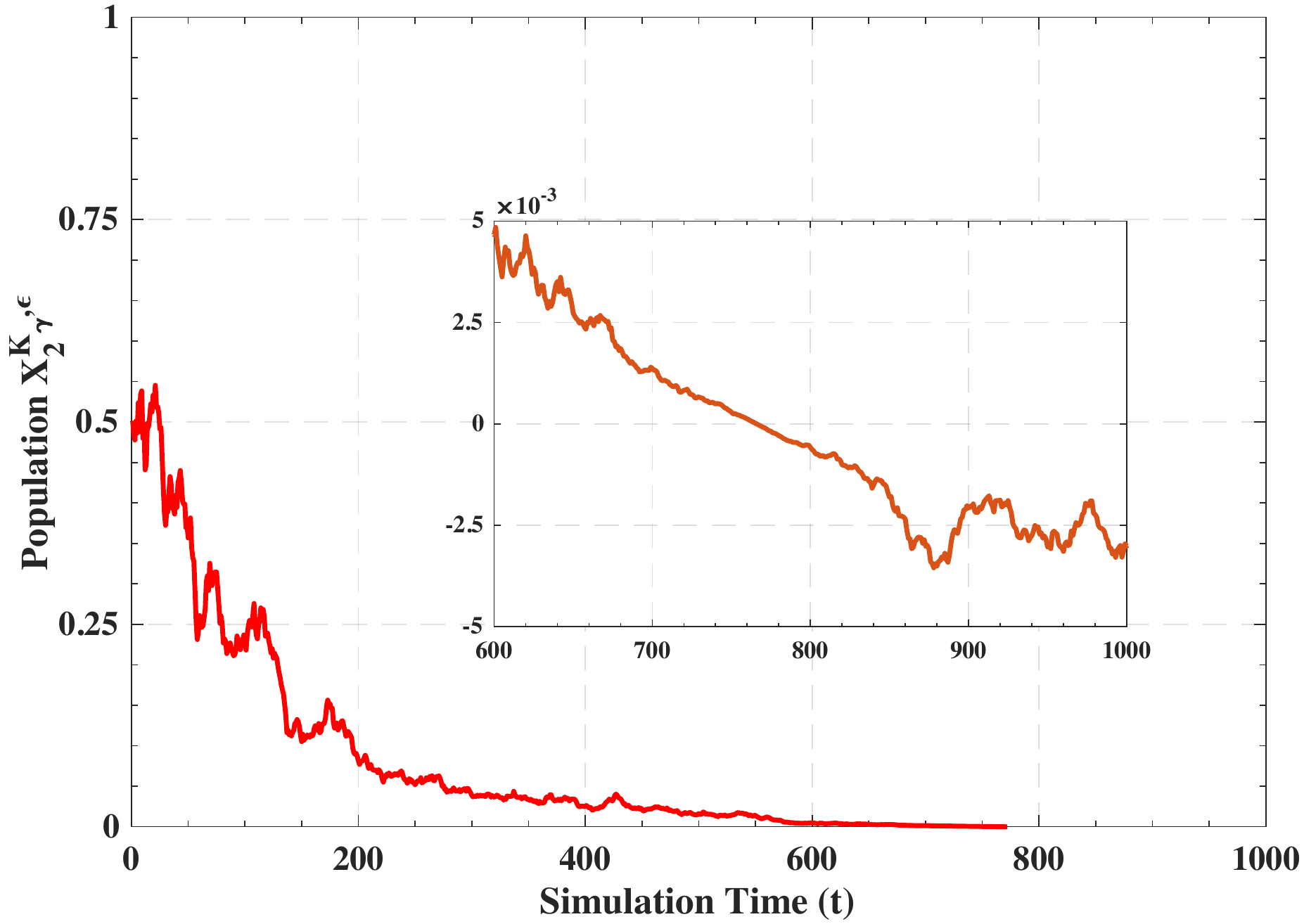}}\\
{\includegraphics[scale=0.311]{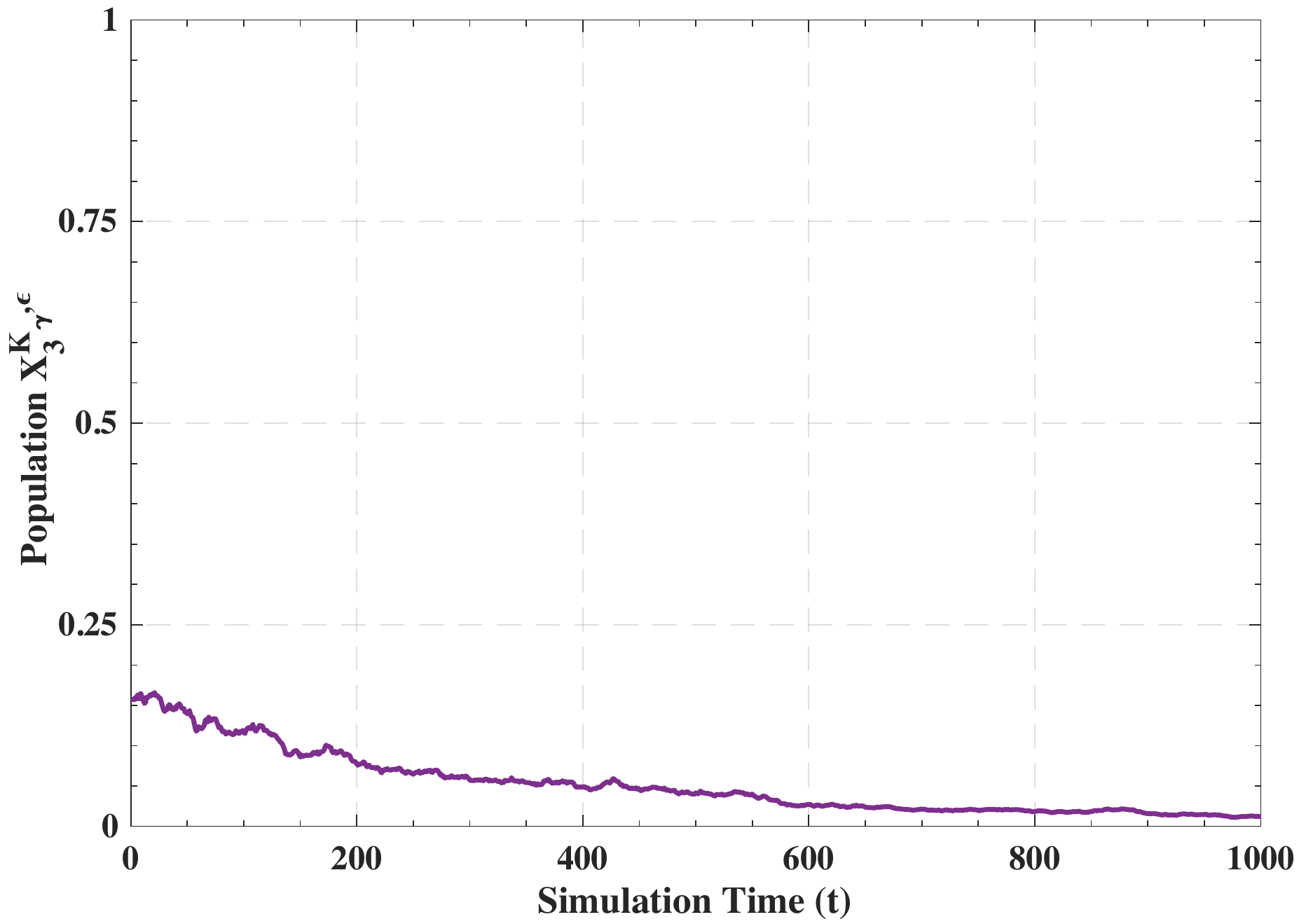}}
{\includegraphics[scale=0.311]{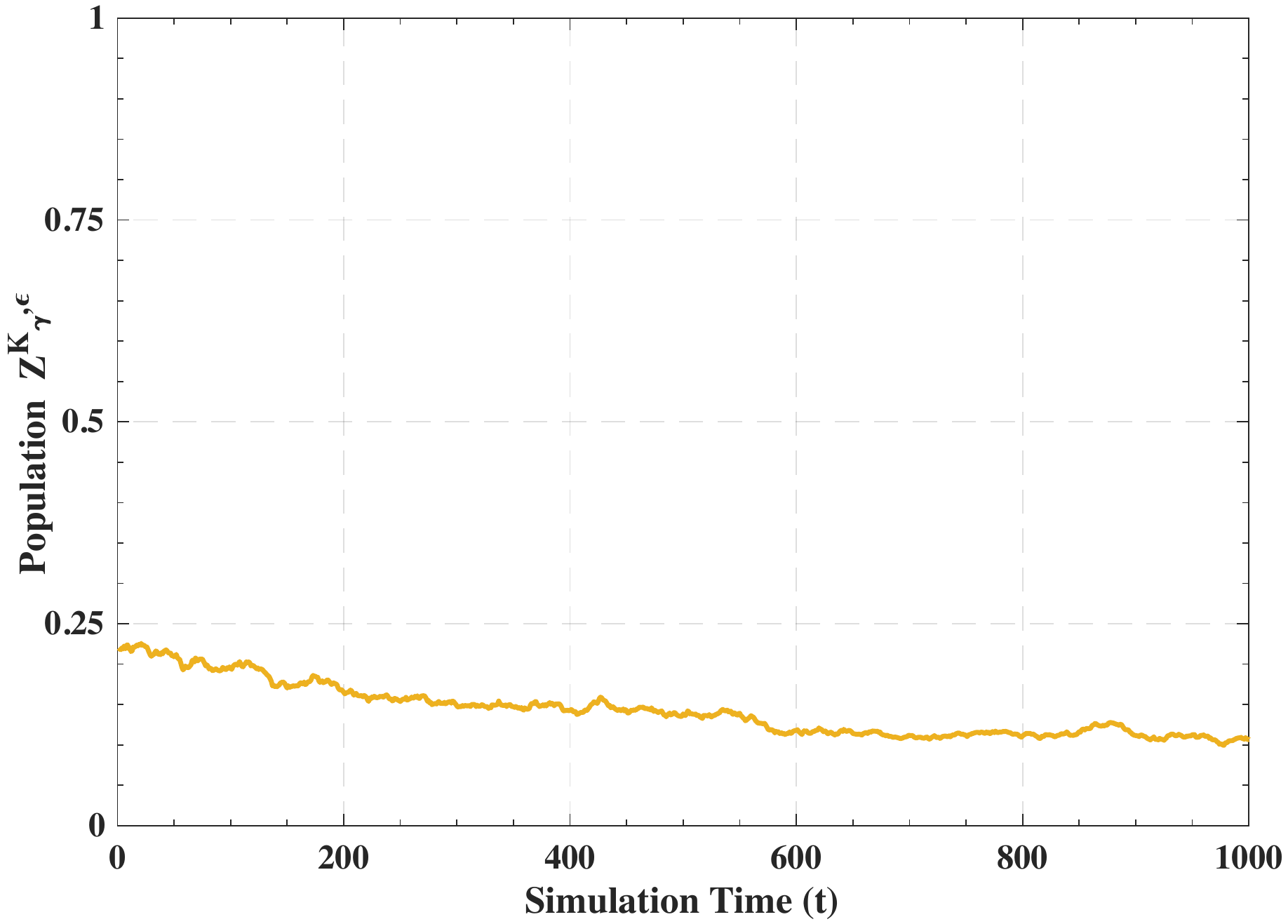}}
\vspace{-1.0mm}
\caption{Population trajectory for small randomly perturbing noise, with an intensity level of $\epsilon = 0.01$.} \label{FG1}
\end{center}
\end{figure}

Then, based on the above discussion, for small random perturbing noise, with an intensity level of $\epsilon = 0.01$, we computed $K_{\tilde{\gamma}}$ using the parameter values from Table~\ref{TB2}, with $\varepsilon=3$ and $\zeta=0.25$, where the resulting controller is given by
\begin{align*}
\tilde{v}(\mathbf{X}(t)) &= \left[ \begin{array}{r r r} -0.1584 & 0.1492 & 0.1422\\
  0.0142 & -2.1721 & -1.9964 \end{array} \right] {\mathbf{X}}(t)\\
    &\triangleq K_{\tilde{\gamma}} {\mathbf{X}}(t), \quad \text{with}\quad \tilde{\gamma} = 0.001.
\end{align*}
Moreover, the corresponding logarithm residence time for the controlled-linear SDE is calculated to be $\tilde{\phi}(D, K_{\tilde{\gamma}}) = 2.4797$. Figure~\ref{FG1} shows the solutions $\bigl({\mathbf{X}}_{0,\mathbf{x}}^{K_{\tilde{\gamma}},\epsilon}(t)\bigr)_{t \ge 0}$ and $\bigl(Z_{0,\mathbf{x}}^{K_{\tilde{\gamma}},\epsilon}(t)\bigr)_{t \ge 0}$ starting from an initial condition $\mathbf{X}(0) = (0.1185,\, 0.5015,\, 0.16)$, with $Z(0) = 0.22$. Notice that, since the addiction-free equilibrium is asymptotically stable, then any solution for the prescription opioid epidemic model, without any random perturbing noise, starting at those points inside the domain $D$ or near to the addiction-free equilibrium point moves closer to it over time and, hence, the trajectory for the unperturbed system will not leave from the domain $D$, for all $t \ge 0$. On the other hand, looking closer at Figure~\ref{FG1}, we observe that for stochastically perturbed prescription opioid epidemic model, the situation is quite different and any solution that starts at those points inside the domain $D$ or near to the addiction-free equilibrium point will eventually leave from the domain $D$. For example, see the figure on the right top side, where the population trajectory $\mathbf{X}_{0,\mathbf{x}}^{K_{\tilde{\gamma}},\epsilon}(t)$ exits from the domain $D$, i.e., $\mathbf{X}_{0,\mathbf{x}}^{K_{\tilde{\gamma}},\epsilon}(t) \notin D$, for $t \ge \tau_{D}^{\epsilon}$ (namely, $X_2^{K_{\tilde{\gamma}},\epsilon}(t) \le 0$, when $t \ge \tau_{D}^{\epsilon}$) with the corresponding simulation exit-time $\tau_{D}^{\epsilon} \approx 781$ (simulation time). 

\section{Concluding remarks} \label{S5}
In this paper, we considered the problem of optimal residence time control for the prescription opioid epidemic dynamical model with small random perturbations. In particular, we argued that such an optimal control problem can be posed as minimizing the asymptotic exit-rate (or equivalently maximizing the expected exit-time) with which the controlled-diffusion process associated with stochastically perturbed epidemic model exits from the given bounded domain and, as a result of this, we established a connection with a controlled-eigenvalue problem. Moreover, we also determined the corresponding admissible optimal Markov control for the HJB equation that has an interpretation of optimally effected removal of susceptible individuals from the population of prescription opioid epidemics or when the addictive users is strategically influenced due to an effective or a more accessible addiction treatment facility. Note that, by considering random processes as perturbations in the prescription opioid epidemic model, we are able to extend the stationarity nature of perturbations and we further provided sufficient information on the probabilistic characteristic of the most probable population trajectory when exiting from the domain of interest, where such additional information could be useful for interpreting outcome-results from opioid-related intervention policies.

\end{document}